\documentclass[10pt]{amsart}
\usepackage{amsmath, amssymb, amsfonts, amscd}
\usepackage{epsfig}

\newcommand{\nc}{\newcommand}
\nc{\vg}{\mathfrak{v} } \nc{\wg}{\mathfrak{w} }
\nc{\zg}{\mathfrak{z} } \nc{\ngo}{\mathfrak{n} }
\nc{\ngoq}{\mathfrak{n}^{\QQ} } \nc{\ngoz}{\mathfrak{n}^{\ZZ} }
\nc{\ggoq}{\mathfrak{g}^{\QQ} } \nc{\kg}{\mathfrak{k} }
\nc{\mg}{\mathfrak{m} } \nc{\bg}{\mathfrak{b} }
\nc{\ggo}{\mathfrak{g} } \nc{\ggob}{\overline{\mathfrak{g}} }
\nc{\sog}{\mathfrak{so} } \nc{\sug}{\mathfrak{su} }
\nc{\spg}{\mathfrak{sp}} \nc{\slg}{\mathfrak{sl} }
\nc{\glg}{\mathfrak{gl} } \nc{\cg}{\mathfrak{c} }
\nc{\rg}{\mathfrak{r} } \nc{\hg}{\mathfrak{h} }
\nc{\tg}{\mathfrak{t} } \nc{\ug}{\mathfrak{u} }
\nc{\dg}{\mathfrak{d} } \nc{\ag}{\mathfrak{a} }
\nc{\pg}{\mathfrak{p} } \nc{\sg}{\mathfrak{s} }
\nc{\lgo}{\mathfrak{l} } \nc{\fg}{\mathfrak{f} }

\nc{\pca}{\mathcal{P}} \nc{\nca}{\mathcal{N}}
\nc{\lca}{\mathcal{L}} \nc{\oca}{\mathcal{O}}
\nc{\mca}{\mathcal{M}} \nc{\tca}{\mathcal{T}}
\nc{\aca}{\mathcal{A}}

\nc{\vp}{\varphi} \nc{\ddt}{\frac{{\rm d}}{{\rm d}t}}
\nc{\im}{\mathtt{i}}

\renewcommand{\l}{\lambda}
\nc{\ala}{Anosov Lie algebra} \nc{\alas}{Anosov Lie algebras}

\nc{\SO}{{\mathrm SO}} \nc{\Spe}{{\mathrm Sp}} \nc{\Sl}{{\mathrm
SL}} \nc{\SU}{{\mathrm SU}} \nc{\Or}{{\mathrm O}} \nc{\U}{{\mathrm
U}} \nc{\Gl}{{\mathrm GL}} \nc{\Se}{{\mathrm S}} \nc{\Cl}{{\mathrm
Cl}} \nc{\Spin}{{\mathrm Spin}} \nc{\Pin}{{\mathrm Pin}}

\nc{\RR}{{\Bbb R}} \nc{\HH}{{\Bbb H}} \nc{\CC}{{\Bbb C}}
\nc{\ZZ}{{\Bbb Z}} \nc{\FF}{{\Bbb F}} \nc{\NN}{{\Bbb N}}
\nc{\QQ}{{\Bbb Q}} \nc{\PP}{{\Bbb P}}

\nc{\vs}{\vspace{.5cm}}
 \nc{\ip}{\langle\cdot,\cdot\rangle}
\nc{\la}{\langle} \nc{\ra}{\rangle} \nc{\unm}{\frac{1}{2}}
\nc{\unc}{\frac{1}{4}} \nc{\und}{\frac{1}{16}} \nc{\f}{\frac}
\nc{\ben}{\begin{enumerate}}\nc{\een}{\end{enumerate}}

 \nc{\no}{\vs\noindent}
\nc{\lam}{\Lambda^2\ggo^*\otimes\ggo} \nc{\tang}{{\rm T}}
\nc{\dif}{{\rm d}} \nc{\preq}{\simeq_K} \nc{\lb}{[\,,\,]}
\nc{\ds}{\displaystyle}
\newcommand{\D}{\mbox{deg}}

\nc{\pf}{\operatorname{pf}} \nc{\ad}{\operatorname{ad}}
\nc{\Ad}{\operatorname{Ad}} \nc{\rank}{\operatorname{rank}}
\nc{\Irr}{\operatorname{Irr}} \nc{\End}{\operatorname{End}}
\nc{\Aut}{\operatorname{Aut}} \nc{\Inn}{\operatorname{Inn}}
\nc{\Der}{\operatorname{Der}} \nc{\Ker}{\operatorname{Ker}}
\nc{\Iso}{\operatorname{I}} \nc{\Diff}{\operatorname{Diff}}
\nc{\Lie}{\operatorname{L}} \nc{\tr}{\operatorname{tr}}

\nc{\dQ}{\operatorname{\mbox{deg}_\mathbb{Q}}}
\nc{\degr}{\operatorname{deg}} \nc{\sen}{\operatorname{sen}}
\nc{\modu}{\operatorname{mod}} \nc{\Ric}{\operatorname{Ric}}
\nc{\sym}{\operatorname{sym}} \nc{\sca}{\operatorname{sc}}
\nc{\scalar}{{\sf s}} \nc{\grad}{\operatorname{grad}}
\nc{\ricci}{\operatorname{ric}} \nc{\Rin}{\operatorname{M}}
\nc{\Le}{\operatorname{L}}
\nc{\level}{\operatorname{level}} \nc{\rad}{\operatorname{r}}
\nc{\abel}{\operatorname{ab}} \nc{\Pf}{\operatorname{Pf}}

\theoremstyle{plain}
\newtheorem{theorem}{Theorem}[section]
\newtheorem{proposition}[theorem]{Proposition}

\newtheorem{lemma}[theorem]{Lemma}
\newtheorem{notation}[theorem]{Notation}

\theoremstyle{definition}
\newtheorem{definition}[theorem]{Definition}

\theoremstyle{remark}
\newtheorem{remark}[theorem]{Remark}

\title[Characterization of 9-dimensional Anosov Lie algebras]{Characterization of 9-dimensional Anosov Lie algebras}

\author{Meera Mainkar,\quad Cynthia E. Will}

\address{Department of Mathematics, Pearce Hall, Central Michigan University, Mt. Pleasant, MI 48858, USA//
 FaMAF and CIEM, Universidad Nacional de C\'ordoba, Haya de la Torre s/n,
5000 C\'ordoba, Argentina } \email{maink1m@cmich.edu,
 cwill@famaf.unc.edu.ar}

\thanks{2000 {\it Mathematics Subject Classification.} Primary: 22E25;
Secondary: 37D20, 20F34. \\
{\it Key words and phrases.}  Anosov Lie algebras,
nilmanifolds, nilpotent Lie algebras,
hyperbolic automorphisms. \\
}

\begin{document}

\maketitle

\begin{abstract}
 The classification of all real and rational Anosov Lie algebras up to dimension $8$  is given by Lauret and Will  \cite{lw2}.
 In this paper we study $9$-dimensional Anosov Lie algebras by using the properties of very special algebraic numbers and Lie algebra classification tools.  We prove that there exists a unique, up to isomorphism, complex $3$-step Anosov Lie algebra of dimension $9$.  In the $2$-step case, we prove that a $2$-step $9$-dimensional Anosov Lie algebra with no abelian factor
 must have a $3$-dimensional derived algebra and we characterize these Lie algebras in terms of their Pfaffian forms.
 Among these Lie algebras, we exhibit a family of infinitely many complex non-isomorphic Anosov Lie algebras.
   \end{abstract}

\section{Introduction}\label{intro}

A diffeomorphism $f$ of a compact differentiable manifold $M$
is called {\it Anosov} if the tangent bundle $\tang M$ admits a
continuous invariant splitting $\tang M=E^+\oplus E^-$ such
that $\dif f$ expands $E^+$ and contracts $E^-$ exponentially. A well-known class of Anosov diffeomorphisms arises as
 follows. Let $N$ be a simply connected nilpotent Lie group and let $\Gamma$ be a discrete subgroup of $N$ such that $N/\Gamma$ is compact. In this case, $N /\Gamma$ is called a \textit{nilmanifold}. If $f$ is a hyperbolic automorphism of $N$ (i.e. no eigenvalue of the differential $d f$ is of absolute value  1) such that
 $f(\Gamma) = \Gamma$, then the induced diffeomorphism $\overline{f}$ on $N / \Gamma$, defined by $\overline{f}(x \Gamma) = f(x) \Gamma$ for all $x \in N$,  is an Anosov diffeomorphism of the nilmanifold $N /\Gamma$. An Anosov diffeomorphism of a nilmanifold $N/\Gamma$ arising in this way is called an \textit{Anosov automorphism} of $N/\Gamma$.
 More generally, one can get examples of Anosov diffeomorphisms on manifolds which are finitely covered by nilmanifolds in the following way.  Let $K$ be a finite group of automorphisms of a simply connected nilpotent Lie group $N$ and let $\Gamma$ be a torsion free discrete subgroup of  $K \ltimes N$ such that  the quotient $N /\Gamma$
  is compact.  Here the action of $\Gamma$ on $N$ is given by $x(\tau, y) =  y \tau(x)$ for $x \in N$ and $(\tau, y ) \in \Gamma$.
   We call the quotient space $N /\Gamma$ an \textit{ infranilmanifold}. If $g$ is a hyperbolic automorphism of $N$ such that $g$ normalizes the subgroup $K$ in the group of automorphisms of $N$ and $g(\Gamma) = \Gamma$, then the induced diffeomorphism on $N/\Gamma$ is called an \textit{ Anosov automorphism of the infranilmanifold $N/\Gamma$. }

In \cite{SS}, S. Smale raised the problem of classifying the
compact manifolds admitting   Anosov diffeomorphisms, and up to
now the only known examples of Anosov diffeomorphisms  are the Anosov automorphisms of infranilmanifolds described above.
 J. Franks \cite{F} and A. Manning \cite{Man} proved that  an Anosov diffeomorphism of a nilmanifold $N/\Gamma$ is topologically conjugate to an Anosov automorphism of   $N/\Gamma$.
This enhances the interest in the problem of
classifying all nilmanifolds which admit Anosov automorphisms.

The first example  of a non-toral nilmanifold
admitting an Anosov automorphism was described by S. Smale
(\cite{SS}). For many years only relatively few examples appeared in the literature, but in recent years
families of nilmanifolds with Anosov automorphisms have been constructed showing that a complete classification does not seem  to be possible (see \cite{L, A-S, Dn, D-M, D, De-Des, lw1, lw2,  M, mw, P}).

Any Anosov automorphism of a nilmanifold $N/\Gamma$ gives rise to a hyperbolic
automorphism $\tau$ of $\ngo$, the Lie algebra of $N$,  
such that  $\tau$ stabilizes a $\mathbb Z$-subalgebra $\Lambda$ of
$\ngo$. Hence the  matrix of $\tau$  with respect to a $\mathbb Z$-basis of $\Lambda$ (i.e. with structure constants in $\ZZ$) lies in
$\Gl_n(\ZZ)$ where $n = \dim \ngo$.

\begin{definition}\label{Anosov} (See \cite{L})
A rational Lie algebra $\ngoq$ (i.e. with structure constants in $\QQ$) of dimension
$n$ is said to be {\it Anosov} if it admits a  hyperbolic automorphism $\tau$
which is {\it unimodular}, i.e.   $[\tau]_{\beta}\in \Gl_n(\ZZ)$ for some basis $\beta$ of $\ngoq$,
where $[\tau]_{\beta}$ denotes the matrix of $\tau$ with respect to the basis $\beta$.  We will  call a
real or complex Lie algebra {\it Anosov} if it admits a rational form which is Anosov.
 \end{definition}

\noindent Therefore, if a nilmanifold $N/\Gamma$ admits an Anosov automorphism then the Lie algebra of $N$ is {\it Anosov.}
 It can be seen that real Anosov Lie algebras give rise to nilmanifolds admitting Anosov automorphisms.
 Hence the problem of classifying nilmanifolds admitting Anosov automorphisms reduces to the classification of real Anosov Lie algebras. So far, in \cite{lw2}, all real and rational Anosov Lie
 algebras of dimension $\leq 8$ are classified up to isomorphism and we note that curiously enough, there are quite a few of them.

\vspace{.5cm}

In this paper we consider the problem of classifying  $9$-dimensional complex Anosov Lie algebras.
 We note that in order to classify  real  Anosov Lie algebras, for each complex Anosov Lie algebras, one has to find all the real forms and study which of them are Anosov.  Therefore, to study complex Anosov Lie algebras is the first step towards the classification.

 In our approach we use two main tools, some algebraic number theory and Lie algebra classification results.
 More precisely, the eigenvalues the hyperbolic automorphism $\tau$ of an Anosov Lie algebra $\ngo$ whose matrix with respect to a $\ZZ$-basis is in $\Gl(n, \ZZ)$, are algebraic units (algebraic integers whose reciprocals are also algebraic integers) of absolute value different from $1$. The properties of these type of algebraic units, studied in \cite{M2}, can be used to give examples or to prove non-existence of Anosov Lie algebras. This approach was introduced in \cite{L} and have also been used in \cite{lw1},\cite{lw2} and \cite{mw}.
 Using this, we prove that every  $9$-dimensional real Anosov Lie algebra, without an abelian factor, is either of type $(6,3)$ or $(3,3,3).$ Recall that

 \begin{definition} \label{type}
{\rm Let $\ngo$ be an $r$-step nilpotent Lie algebra, i.e., the
lower central series $\{C^i(\ngo)\}$ (defined by $C^0(\ngo)=\ngo$ and
$C^i(\ngo)=[\ngo,C^{i-1}(\ngo)]$ for $i \geq 1$) satisfies
$C^{r-1}(\ngo)\not=0$ and $C^r(\ngo)=0$. Then the {\em type}
 of $\ngo$ is the  $r$-tuple of positive
integers $(n_1,\cdots,n_r)$, where $n_i=  \dim C^{i-1}(\ngo) /
C^{i}(\ngo)$.}
\end{definition}

\begin{definition} \label{abelian}
{\rm  Let $\ngo$ be a  Lie algebra. An {\em abelian factor} of
$\ngo$ is an abelian (Lie) ideal $\ag$ of $\ngo$ such that $\ngo = \mg
\oplus \ag$ for some ideal $\mg$ of $\ngo$. }
\end{definition}

We then show that there exist a complex Anosov Lie algebra of type $(3, 3, 3)$ that is unique up to Lie algebra isomorphism. It can be seen  that if $\ngo$ is an Anosov Lie algebra of type $(3, n_2, n_3)$, then $n_2 = 3$  and $n_3 \geq 3$ (see \cite[Proposition 2.3]{lw1}). The type $(3, 3, 3)$ example mentioned above is the first example, to our knowledge, of an Anosov Lie algebra of type $(3,3,*)$ showing that the condition
  $n_3 \ge 3$ is in fact attained. To show that this algebra is actually an Anosov Lie algebra, we use some special algebraic units constructed from the roots of unity in a way that seems to be generalizable.

 The $(6,3)$ case is more involved. By using the classification given by Galitzki and Timashev in \cite{GT},  we show that if a Lie algebra of type  $(6, 3)$ is Anosov then its Pfaffian form is projectively equivalent to $xyz$ or $0$ (see \cite{L}) then we give a list of possible candidates for $(6, 3)$ type complex Anosov Lie algebras (Theorem \ref{teo}). We note that among these algebras there is an infinite family of complex non-isomorphic Anosov Lie algebras (see Proposition \ref{sU_2+U_3anosov} and Remark \ref{infinite}).  It was proved in \cite{lw2} that up to isomorphism there are only finitely many real Anosov Lie algebras up to dimension $8$ (see \cite[Table 3]{lw2}). This shows that the smallest dimension in which there are infinitely many complex non-isomorphic Anosov Lie algebras is $9$.

\section{Preliminaries}

In this section, in order to introduce the framework in which we are going to work, we begin by recalling the following proposition, proved in \cite{lw1}.

\begin{proposition} \label{deco}
 Let $\ngo$ be a real $r$-step nilpotent Lie algebra. We define $C^i(\ngo)$
 inductively by $C^i(\ngo) = [\ngo, C^{i-1}(\ngo)]$ for $i \geq
 1$, where $C^0(\ngo) = \ngo$. If $\ngo$ is an Anosov Lie algebra then
 there exist a decomposition
$\ngo=\ngo_1\oplus \cdots \oplus\ngo_r$ satisfying
$C^i(\ngo)=\ngo_{i+1}\oplus...\oplus\ngo_r$, $i=0,...,r$, and a
hyperbolic $\tau\in\Aut(\ngo)$ such that
\begin{itemize}
\item[(i)] $\tau(\ngo_i)=\ngo_i$ for all $i=1,...,r$.

\item[(ii)] $\tau$ is semisimple (in particular $\tau$ is diagonalizable
over $\CC$).

\item[(iii)] For each $i$, there exists a basis $\beta_i$ of
$\ngo_i$ such that $[\tau_i]_{\beta_i}\in\Sl_{n_i}(\ZZ),$ where
$n_i=\dim{\ngo_i}$ and $\tau_i=\tau|_{\ngo_i}$.
\end{itemize}
\end{proposition}

Let  $\ngo$ and $\tau$  be as in Proposition \ref{deco},
 i.e. $\tau$  is  a hyperbolic automorphism of $\ngo$ and there exists a $\ZZ$-basis
  of $\ngo,$ $\beta,$ such that $[\tau]_{\beta} \in  \Sl_n(\ZZ)$ where $n$ = dim $\ngo$.  We note that  the eigenvalues of $\tau$ are algebraic units, that is
 each eigenvalue of $\tau$ satisfies
 a monic polynomial equation with integer coefficients and constant term $1$. This follows from the existence of a basis
 of $\ngo$ with respect to which the matrix of $\tau$ is in $ \Sl_{n}(\ZZ).$

 For an algebraic number $\lambda \in  \mathbb C$, we denote by $\D(\lambda)$
 the degree of $m_\l(x)$, the  irreducible monic polynomial over $\QQ$ annihilated by
 $\lambda$  and by the
 {\it  conjugate of $\l$} 
 the other roots of $m_\l(x).$  We then have the following lemma:

\begin{lemma}\label{util}
Let $\ngo$ be a real  $r$-step Anosov  Lie algebra, and
let $\tau$ and $\ngo=\ngo_1\oplus\ngo_2\oplus\dots \oplus\ngo_r$  be
as in {\rm Proposition \ref{deco}}.  Let  $n_i=\dim{\ngo_i}$ and $\tau_i=\tau|_{\ngo_i}$ for $1 \leq i \leq r$.
  Then every eigenvalue  $\lambda_i$   of $\tau_i$ is an  algebraic
unit and  $1< \D(\l_i) \le n_i$ for all $i$.

\end{lemma}

The following definition will be used in  the following sections.

\begin{definition}\label{splitting} Let $V$ be a real vector space of dimension $n$. Let $\sigma$ be a linear automorphism of $V$ whose characteristic
polynomial has integer coefficients.  We say that  $\sigma$ has a {\it splitting }
 $[k_1; \dots; k_m],$ where $k_i \in \NN, k_i \geq k_{i+1},$ if the characteristic polynomial of $\sigma$ can be
 written as a product of $m$ irreducible polynomials (over $\ZZ$)
 $f_1, f_2, \ldots f_m$ such that $\degr f_i = k_i$ for all $i$.
\end{definition}

For example if $\sigma: \mathbb R^3 \rightarrow \mathbb R^3$ is a linear automorphism given by $\sigma(x, y, z) = (2x+ y, x+y, 3z)$ for all
$x, y, z \in \mathbb R$, then the characteristic polynomial of $\sigma$ is $X^3-6X^2+10X-3$ with irreducible factors over $\ZZ$,  $X^2-3X+1$ and $X-3$. Hence the splitting of $\sigma$ is $[2;1]$.

In our setup, if $\ngo, \tau, \ngo_i$ and $\tau_i$ are as in Proposition \ref{deco}, we are going to look at the possible
splittings  of $\tau_i'$s. Note that for a fixed $i$,  the characteristic polynomial of $\tau_i$ has integer coefficients (Proposition
 \ref{deco} (iii)). If $\tau_i$ has a splitting  $[k_1; \dots; k_m]$, then $k_1+ \dots + k_m=n_i$  where $n_i = \dim{\ngo_i}$.
  Since $\tau_i$ is hyperbolic and $[\tau_i]_{\beta_i}\in\Sl_{n_i}(\ZZ),$ for some basis $\beta_i$ of $\ngo_i$, we also have that $k_j \ne 1 $ for all $j$ (see also \cite{lw1} Appendix).

We note that if $i >1$ each eigenvalue of $\tau_i$ is a product of certain eigenvalues of $\tau_j$ with $j<i$ and moreover, it is an algebraic unit.

We will investigate how the product of such special algebraic numbers can behave.
It is proved in \cite{M2} that this behavior can not be too wild. In fact, one has the following lemma that can be deduced from
\cite{M2}.

\begin{lemma}\label{m1}
Let $\tau$ be an Anosov automorphism of a nilpotent Lie algebra $\ngo$.  Let $\alpha$ and $\beta$ be
two eigenvalues of $\tau$. Then the following hold:
 \ben \item[1.] If $\D(\alpha)$ and
 $\D(\beta)$ are relatively prime then  $\D(\alpha \beta)$ cannot be a prime.
\item[2.] If $\alpha\beta$ is also an eigenvalue of $\tau$ and $\D
(\alpha) \le \D(\beta)$ then $$\mbox{g.c.d.}(\D(\beta), \D(\alpha\beta))\ne 1.$$
\een\end{lemma}

\proof If $\D(\alpha)$ and
 $\D(\beta)$ are relatively prime and   $\D(\alpha \beta)$ is prime, then it follows from Corollary 2 of \cite{M2} that either $|\alpha| = 1$ or $|\beta| = 1$. This contradicts  our assumption that both $\alpha$ and $\beta$ are eigenvalues of a hyperbolic automorphism $\tau$.

Suppose that $\alpha \beta$ is an eigenvalue of $\tau$, $\D(\alpha) = \D(\beta)$ and  $\mbox{g.c.d.}(\D(\beta), \D(\alpha\beta))= 1$.  Then $|\alpha \beta| = 1$ by Corollary 1 of \cite{M2} which is a contradiction. Also if  $\D
(\alpha) < \D(\beta)$, then  by Lemma 2 of \cite{M2} we have  $\mbox{g.c.d.}(\D(\beta), \D(\alpha\beta))\ne 1.$ \hfill$\square$
\vspace{.5cm}

We recall  \cite[Theorem 3.1]{lw1}: Let $\ngo$ be a rational Lie
algebra and let $\ngo = \widetilde{\ngo} \oplus  \mg$ be a Lie direct sum,
where $\mg$ is a  maximal abelian factor
  of $\ngo$. Then $\ngo$ is Anosov if and only if
$\widetilde{\ngo}$ is Anosov and dim $\mg \geq 2$. In view of this, we are
interested in studying  Anosov Lie algebras without an abelian
factor (see Definition \ref{abelian}).

For an  Anosov Lie algebra $\ngo$ without an abelian factor, we will study the  properties of  eigenvalues of an Anosov automorphism of $\ngo$  and by using these properties we will deduce the Lie bracket structure on
$\ngo_{\CC} = \ngo\otimes\CC$. Moreover, we  also  note that if $\ngo_{\CC}$ has an abelian factor,
then $\ngo$ must also carry an abelian
factor. In fact, $\zg(\ngo_{\CC}) \cap [\ngo_{\CC}, \ngo_{\CC}] =
 (\zg(\ngo) \cap [\ngo, \ngo])_{\CC}$ and hence if $\zg(\ngo_{\CC}) \cap [\ngo_{\CC}, \ngo_{\CC}]
  \neq \zg(\ngo_{\CC})$ then $(\zg(\ngo) \cap [\ngo, \ngo]) \neq
  \zg(\ngo)$ where $\zg(\ngo)$ denotes the center of $\ngo$.

\vspace{.5cm}

Finally, we state the following lemma that can be proved by using Proposition \ref{deco}.

\begin{lemma}{\label{rest}}Let $\ngo$ be an Anosov Lie algebra
of type $(n_1, \dots,n_r)$ and let
$\ngo=\ngo_1\oplus...\oplus\ngo_r$ be the decomposition of $\ngo$
given in Proposition \ref{deco}. Then $\widetilde{\ngo}=\ngo / \ngo_r$  is Anosov.
\end{lemma}

\proof If  $\ngo$ is an Anosov Lie algebra, $\tau$ and $\ngo=\ngo_1\oplus...\oplus\ngo_r$ are as in Proposition \ref{deco}
then,
 since $\tau\ngo_i=\ngo_i$ for all $i=1,...,r$, it is easy to see that it induces an automorphism of $\widetilde{\ngo}=\ngo / \ngo_r \simeq \ngo_1\oplus...\oplus\ngo_{r-1}.$
 Also this automorphism is hyperbolic and to see that it is   unimodular, recall that for each $i$, there exists a basis $\beta_i$ of
$\ngo_i$ such that $[\tau_i]_{\beta_i}\in\Sl_{n_i}(\ZZ),$ where
$n_i=\dim{\ngo_i}$ and $\tau_i=\tau|_{\ngo_i}$. \hfill$\square$

\vspace{.3cm}

Note that this argument is valid not only at the real or complex level but also in the rational case.
  We also note that if
$\ngo$ is a three step nilpotent Anosov Lie algebra,
 then $\widetilde{\ngo}$ is 2-step and  the decomposition $\widetilde{\ngo}=\ngo / \ngo_3 \simeq \ngo_1\oplus \ngo_2$ gives the type of  $\widetilde{\ngo}.$

\vspace{.5cm}

\section{Dimension 9}

\vspace{.3cm}

In this section we will study  Anosov Lie algebras of
dimension $9.$
 We will prove that an Anosov Lie algebra without an abelian factor  must  be of type $(6,3)$ or $(3, 3, 3)$.  Moreover,  we will prove that there is only one complex Anosov Lie algebra (up to a Lie algebra isomorphism)
without an abelian factor of type $(3, 3, 3)$. The case $(6, 3)$  is the hardest.

Using  \cite[Proposition 2.3]{lw1} we can see that the
possible types for a $9$-dimensional Anosov Lie algebra are
$(7,2),\,\mathrm{(6,3)},\,(5,4),\,(4,5),\,(5,2,2),\,(4,3,2),\,(4,2,3)$
and $\mathrm{(3,3,3)}$.
As a corollary of Lemma \ref{rest} and the fact that there is no
non-toral $7$-dimensional Anosov Lie algebra (see \cite{lw2}), we get that there
are no Anosov Lie algebras of type $(5, 2, 2)$ and $(4, 3, 2)$.
Finally, it is not hard to see by using Lemma \ref{m1}, that there are no Anosov Lie algebras with no abelian factor of type $(7,2),$ $(5,4),$ $(4,5)$ and $(4,2,3).$ We will illustrate by looking at the case $(5,4)$.

In the rest of this section, we will use the following notation:

\begin{notation}\label{eigenvectors}{\rm  If $\ngo$ is an Anosov Lie algebra and $\tau,$ $\tau_i$ and $\ngo_i$
 are  as in Proposition \ref{deco},
 we will denote the eigenvalues of $\tau_1$ by $\lambda_i'$s, the eigenvalues of $\tau_2$ by $\mu_j'$s and the eigenvalues of $\tau_3$  by $\nu_k'$s, and the corresponding eigenvectors by $X_i'$s,
  $Y_j'$s and  $Z_k'$s  (in  $\ngo_{\CC}$) respectively.}
\end{notation}
\vspace{.2cm}

\begin{proposition} \label{(5,4)}
 There is no Anosov Lie algebra of type $(5,4)$.
\end{proposition}

\proof Suppose that there exists
 an Anosov Lie algebra $\ngo$ of  type $(5, 4)$, and let $\tau$ be an Anosov automorphism as in Proposition \ref{deco}. We note that the possible splittings
(see Definition \ref{splitting}) of $\tau_1$ are $[5]$ and $[3;2].$ In the first case, for each nonzero
bracket among the eigenvectors of $\tau_1$ we get an eigenvector of $\tau_2$ (of degree $2$ or $4$). That is if $[X_i,X_j] \ne 0$ for some $1 \le i,j \le 5,$ we then have that $\l_i \l_j=\mu_k$ has degree $2$ or $4.$ Either way this contradicts Lemma \ref{m1} since $g.c.d(5,4)=g.c.d(5,2)=1.$
On the other hand, if the splitting is $[3;2]$  by the same argument, we get that
$[X_i,X_j]=0$ for $i,j \in \{1,2,3\}$, and if $\mu_k=\l_i \l_j$
with $i \le 3$ and $j=4$ or $5,$
 then g.c.d.$(\degr\l_i, \degr\mu_k)=1$, contradicting again Lemma \ref{m1}. This means that $[X_i,X_j]=0$ for all $i, j$ since $\l_4 = \l_5^{-1}$ and therefore $[X_4,X_5]=0.$ This is a contradiction because $\ngo$ is of type $(5, 4)$. \hfill$\square$


\vspace{.3cm}

\begin{proposition} \label{(3,3,3)}
 There exists only one complex  Anosov Lie algebra of type $(3, 3, 3)$ up to Lie algebra isomorphism.
\end{proposition}

\proof  Let $\ngo$ be a nilpotent Lie algebra of type $(3,3,3)$ that admits an Anosov automorphism $\tau$ as in Proposition \ref{deco}. In this case we note that the characteristic polynomials of $\tau_1$, $\tau_2$ and $\tau_3$ are all irreducible degree 3 polynomials over $\QQ$.
Following Notation \ref{eigenvectors}, we can assume with no loss of generality, that in $\ngo_{\CC}$ (complexification of $\ngo$):
$$\begin{array}{lcr}
 [X_1,X_2]=Y_1,& [X_2,X_3]=Y_2,& [X_1,X_3]=Y_3,
 \end{array}
 $$
and since $\l_1 \l_2 \l_3 = \pm 1$ we also have that
\begin{equation}\label{3331}
\begin{array}{lcr}
 [X_3,Y_1]=0,& [X_1,Y_2]=0, & [X_2,Y_3]=0.
 \end{array}
  \end{equation}

 Let $\ngo(a,b,c)$ denotes the Lie algebra given by
 \begin{equation}\label{n333}
\begin{array}{lcr}
 [X_1,X_2]=Y_1,& [X_2,X_3]=Y_2,& [X_1,X_3]=Y_3,\\
& & \\

 [X_1,Y_1]=a\;Z_1,& [X_2,Y_2]=b\;Z_2,& [X_3,Y_3]=c\;Z_3,

 \end{array}
 \end{equation}

 for $a, b, c \in \CC.$  If $abc \neq 0$,  it can be seen that $\ngo(a,b,c)$ is isomorphic to $\ngo(1,1,1)$ by changing the basis in the center.
  To see that $\ngo_{\CC}$ is isomorphic to $\ngo(1,1,1),$ note that otherwise we have (in $\ngo_{\CC}$)
\begin{equation}\label{3332}
[X_i,Y_j]=a\,Z_k, \qquad \mbox{and}\qquad  [X_i,Y_l]=b\,Z_r
\end{equation}
 for some non zero $a,b\in \CC$ and some $i$ that we can assume to be $1$.
It is clear from (\ref{3331}) that $j,l \ne 2$ so we may assume
that $j=1$ and $l=3.$  Hence, if $k=r,$ (\ref{3332}) implies that
$\l_1^2\l_2=\l_1^2\l_3$ and then
$\l_2=\l_3.$  This is a contradiction because $\lambda_i'$s are roots of the characteristic polynomial of $\tau_1$ and they have
to be distinct because the characteristic polynomial of $\tau_1$ is irreducible over $\QQ$.
On the other hand, if $k \ne r,$ with no loss of
generality we can assume that
$$\begin{array}{lr}
 [X_1,Y_1]=a\;Z_1,& [X_1,Y_3]=b\;Z_2.
 \end{array}$$
Since $Z_3$ is an eigenvector of $\tau_3$ corresponding to $\nu_3$ (see Notation \ref{eigenvectors}),
 $Z_3=[X_i,Y_j]$ for some $i,j$. Using (\ref{3331}) and that $\nu_k'$s are all distinct, the  possibilities for $(i,j)$ are $(2,1),\;(2,2),\;(3,2),$
 and $(3,3)$ but all of them lead  us to the same kind of contradiction by checking that the product of the eigenvalues
 in the center equals $\nu_1^2$, $\nu_1$, $\nu_2$, $\nu_2^2$ respectively. Indeed, if for example we consider case
 $(2,2),$ we obtain
  $$\begin{array}{lcr}
 [X_1,Y_1]=a\;Z_1,& [X_1,Y_3]=b\;Z_2,& [X_2,Y_2]= \;Z_3,
 \end{array}$$
 \noindent and therefore $1=\l_1^2\l_2 . \l_1^2\l_3 . \l_2^2\l_3 = \l_1^2\l_2 = \nu_1, $
contradicting the fact that $\tau_3$ is hyperbolic. This shows that $\ngo_{\CC}$ is isomorphic to $\ngo(1,1,1)$.

In the following we will show that $\ngo(1,1,1)$ is Anosov by constructing a hyperbolic automorphism $\sigma$ and a $\ZZ$-basis of $\ngo(1,1,1)$ preserved by $\sigma$ such that the matrix of $\sigma$ in that basis has integer entries.

Consider the polynomial in $\ZZ[X]$ given by
\begin{equation}\label{degree3}
  f(X) = X^3 - 3X + 1.
\end{equation}
Then its roots are given by
  \begin{equation}\label{lambdas333}
  \lambda_1=\xi + \xi^8, \quad \lambda_2=\xi^2 + \xi^7,\quad \lambda_3=\xi^4 + \xi^5,
  \end{equation}
 \noindent where $\xi=e^{2i\pi/9}$ a ninth root of unity.
 In this case we have that the extension $\QQ(\lambda_1)$ is a cyclic extension of degree $3$
 over $\QQ$ (see \cite[p. 543]{A}) and moreover, a straightforward calculation shows that
 \begin{equation}\label{lambdascuad}
  \lambda_1=\l_3^2-2 \quad \lambda_2=\l_1^2-2,\quad \lambda_3=\l_2^2-2.
  \end{equation}

Let
$$ \begin{array}{lll}
\mu_1 = \l_1\l_2,& \mu_2 = \l_2\l_3,& \mu_3 = \l_1\l_3,\\
&&\\
\nu_1=\l_1 \mu_1,& \nu_2=\l_2 \mu_2,& \nu_3=\l_3 \mu_3.
\end{array}$$
We consider now the automorphism $\sigma$ as above, corresponding to
these $\l_i$'s, that is, defined by $\sigma(X_i) = \l_i X_i,\; \sigma(Y_i) = \mu_i Y_i$ and $\sigma(Z_i) = \nu_i Z_i$
 for $i = 1, 2, 3$.

  We note that $\mu_i, \nu_i \notin \QQ$  for all $1 \le i \le 3$.  For, if $\mu_1 = \lambda_1\lambda_2 \in \QQ$, then $\mu_1  = \pm 1$ because $\mu_1$ is an algebraic unit. Then we get a contradiction that $\lambda_3 = \pm 1$ since $\lambda_1\lambda_2\lambda_3 = 1$. Similarly if $\nu_1 \in \QQ$, then $\nu = \pm 1$ and $\lambda_1 = \pm \lambda_3$ which is a contradiction.  Now since $\mu_i, \nu_i  \in \QQ(\lambda_1)$ and $\QQ(\lambda_1)$ is an extension of degree $3$ over $\QQ$,  $\D(\mu_i) = 3 = \D(\nu_i)$ for all $1 \le i \le 3$.

Concerning the new basis, for $i=1, 2, 3$ let us denote by
$$ \begin{array}{l}
\mathcal{X}_i = \displaystyle{\sum_{j= 1}^{3}} \l_j^{i-1} \, X_j, \\
\mathcal{Y}_i = \l_3^{1 - i}(\l_2-\l_1)Y_1+ \l_1^{1-i}(\l_3-\l_2)Y_2+\l_2^{1-i}(\l_3-\l_1)Y_3,\qquad  \\ \\
\mathcal{Z}_i = (\l_1)^{i-1} (\l_2-\l_1)Z_1+ (\l_2)^{i-1}(\l_3-\l_2)Z_2+ (\l_3)^{i-1}(\l_3-\l_1)Z_3.
\end{array}
$$

Let
 $\beta_1=\{\mathcal{X}_i, \,1\le i \le 3\},$
$\beta_2=\{\mathcal{Y}_i, \,1\le i \le 3 \}$ and
$\beta_3=\{\mathcal{Z}_i, \,1\le i \le 3 \}$. It is easy to see that
$\beta_i'$s are
 linearly independent sets.   Moreover, we can see that
$\sigma(\mathcal{X}_i)=\mathcal{X}_{i+1}$,
$\sigma(\mathcal{Y}_i)=\mathcal{Y}_{i+1}$ for $i=1,2$, $\sigma(\mathcal{X}_3)
=
 3 \mathcal{X}_2 - \mathcal{X}_1$ and $\sigma(\mathcal{Y}_3) = 3
\mathcal{Y}_3 - \mathcal{Y}_1$ by  using that $f(\l_i)=0$ and
therefore $\l_i^{-1}=3-\l_i^2$. Also, by using (\ref{lambdascuad})
one can  see that $\sigma(\mathcal{Z}_i)$ is an integer linear
combination of the $\mathcal{Z}_i's$ for all $i$.  In fact, for
example
$$\begin{array}{l}
\sigma(\mathcal{Z}_1)= \l_1^2 \l_2(\l_2-\l_1)Z_1+ \l_2 ^2\l_3(\l_3-\l_2)Z_2+ \l_3^2\l_1(\l_3-\l_1)Z_3
 \\ \\
\quad =\l_1 ^2 (\l_1^2-2)(\l_2-\l_1)Z_1+ \l_2^2(\l_2^2-2)(\l_3-\l_2)Z_2+ \l_3^2(\l_3^2-2)(\l_3-\l_1)Z_3\\ \\
\quad=(\l_1 ^2 -\l_1)(\l_2-\l_1)Z_1+ (\l_2 ^2 - \l_2)(\l_3-\l_2)Z_2+ (\l_3^2 - \l_3)(\l_3-\l_1)Z_3 \\ \\
\quad= \mathcal{Z}_3-\mathcal{Z}_2.
\end{array}$$

Therefore, $\beta = \{
\mathcal{X}_i, \mathcal{Y}_j, \mathcal{Z}_k: i, j, k = 1, 2, 3\}$ is a basis of $\ngo(1,1,1)$ preserved by $\sigma$ and moreover $[\sigma]_\beta \in GL(9, \ZZ)$.
 To conclude, it remains to show that $\beta$ is also a $\ZZ$-basis.
We first note that since $\l_1\l_2\l_3 = 1$,  $\l_1 + \l_2 \l_3 = 0$
and  $\l_i^{-2} = 3 \l_i^{-1} - \l_i$, we have that
\begin{equation}\label{eq3333}
\begin{array}{ll} [\mathcal{X}_1,\mathcal{X}_2]=\mathcal{Y}_1, &
[\mathcal{X}_1,\mathcal{X}_3]=\mathcal{Y}_3-3\mathcal{Y}_2.
\end{array}
\end{equation}
  Moreover, by using (\ref{lambdascuad}) as above, we have
$$\begin{array}{lcl}
 [\mathcal{X}_1,\mathcal{Y}_1] =\mathcal{Z}_1, &
 [\mathcal{X}_1,\mathcal{Y}_2] =\mathcal{Z}_2 - \mathcal{Z}_1, &
[\mathcal{X}_1,\mathcal{Y}_3] = \mathcal{Z}_3 - 2\mathcal{Z}_2 +
\mathcal{Z}_1.
\end{array}$$

Note that since $\sigma$ in an automorphism of $\ngo$ and by our
construction of the basis, these  are the only Lie brackets one needs to
check. Hence $\beta$ is a $\ZZ$-basis of $\ngo(1,1,1)$ preserved by
the hyperbolic automorphism $\sigma,$ such that $[\sigma]_\beta \in GL(9,
\ZZ)$, and therefore $\ngo(1,1,1)$ is an Anosov Lie algebra. \hfill$\square$
\vspace{0.5cm}

\begin{remark} {\rm We will use the above  proof as a
model proof to show that a given Lie algebra is Anosov in some of
the  following cases  where we will  just exhibit the automorphism
and the $\ZZ$-basis.}
\end{remark}

\begin{remark} {\rm Note that with this result we have shown that there is only one $k$-step nilpotent
complex Anosov Lie algebra  of dimension 9 with $k > 2$ (see the beginning of this section). }
\end{remark}


\vspace{.3cm}


\noindent {\bf Type $(6, 3)$}. To study Anosov Lie algebras of type $(6, 3)$,  we will use the classification of complex nilpotent Lie algebras of type $(6, 3)$ by Galitzki and Timashev given in \cite{GT}. We will begin by dividing the study according to the possible splittings of the Anosov automorphism to get some condition on the bracket structure. Then, using classification from \cite{GT},
  we give a list of  non-isomorphic Lie algebras of type $(6,3)$ as possible candidates for Anosov algebras.
  In most cases we will give a  rational basis to show that they are indeed Anosov Lie algebras.

\begin{notation}\label{Unotation}
{\rm  To refer to the classification given in \cite{GT} we will introduce here notation for some Lie
algebras from \cite[Section 4]{GT}  which will be used in this section. For $s_1,\,s_2,\,s_3 \in \CC$, we will denote by
$s_1\;\mathcal{U}_1 + s_2\;\mathcal{U}_2 + s_3\;\mathcal{U}_3$ a  2-step nilpotent Lie algebra of type $(6, 3)$ with a basis $\{ v_i : 1 \le i \le 6 \} \cup \{w_j : 1 \le j \le 3\}$ and nonzero Lie brackets given by

\begin{equation}\label{63rst}
 \begin{array}{lcl}
 [v_1,v_2]=s_1\; w_1 & [v_3,v_4]=s_1\;w_2 & [v_5,v_6]=s_1\;w_3 \\

 [v_5,v_4]=s_2\; w_1 & [v_1,v_6]=s_2\; w_2 & [v_3, v_2]=s_2\; w_3 \\

 [v_3,v_6]=s_3\; w_1 & [v_5,v_2]=s_3\;w_2 & [v_1,v_4]=s_3\;w_3.
 \end{array}
 \end{equation}

 Here, $\mathcal{U}_i$ is the $9$-dimensional nilpotent Lie algebra  with Lie bracket given by the $i^{th}$ row of
(\ref{63rst}) and with $s_i = 1$.

Let $V$ be a $6$-dimensional vector space and let $W$ be a $3$-dimensional vector space. The classification of nilpotent Lie algebras of type $(6, 3)$
has been given in \cite{GT} by classifying the tensors in $\wedge^2V  \otimes W$ under the action of $SL(V) \times SL(W)$. Each tensor is decomposed into a semisimple and a nilpotent part. 
They obtained $7$ families that are divided according to the semisimple part which is given in terms of $s_1\mathcal{U}_1
 + s_2 \mathcal{U}_2 + s_3 \mathcal{U}_3$ for some $s_1, s_2, s_3.$ The nilpotent parts are listed in tables which we are going to denote by
   $\mg_j$ where $j$ is its number in the corresponding table.  For example, by $\mathcal{U}_1 + \mg_{11}$ from Family $6$ (and Table $7$) from
\cite[Section 4]{GT}, we denote the $9$-dimensional Lie algebra with a basis  $\{ v_i, w_j: 1 \le i \le 6, 1 \le j \le 3\}$  and nonzero Lie brackets given by

 \begin{equation}\label{U_1+m_11}
 \begin{array}{lcl}
 [v_1,v_2]=\; w_1 &  [v_3, v_4]=\;w_2  & [v_5,v_6]=\;w_3 \\

 [v_1,v_5]=\; w_2  & [v_3,v_6]=\; w_1

 \end{array}
 \end{equation}

Note that since we are
interested in Lie algebras of type $(6,3)$ up to isomorphism, we need orbits of tensors
under  the action of $\Gl(V)\times\Gl(W),$ and therefore to get this classification, the canonical form for the semisimple part of the tensors can be reduced by multiplying by a nonzero scalar. Then, for example we have that
$s\,\mathcal{U}_2+r\,\mathcal{U}_3 \simeq s'\,\mathcal{U}_2+\mathcal{U}_3$ for any $r \ne 0$ and the
semisimple part of Family $4$ (and $5$) are all isomorphic (see \cite{GT}).

Next, for completeness, we will list the $\mg_{j}$'s we are actually using in this section, by specifying the nonzero Lie brackets of their basis vectors  $\{ v_i, w_j: 1 \le i \le 6, 1 \le j \le 3\}.$

 From Family $6$ and Table $7$ :

\begin{equation}\label{m_6}\tag{$\mg_{6}$}
 \begin{array}{lcl}
 [v_1,v_4]=\; w_3  & [v_1,v_6]=\; w_2 & [v_3, v_5] =\, w_1

 \end{array}
 \end{equation}

\begin{equation}\label{m_7}\tag{$\mg_{7}$}
 \begin{array}{lcl}
 [v_1,v_3]=\; w_3  & [v_1,v_5]=\; w_2 & [v_3, v_6] =\, w_1

 \end{array}
 \end{equation}

\begin{equation}\label{m_10}\tag{$\mg_{10}$}
 \begin{array}{lcl}

 [v_1,v_3]=\; w_3  & [v_1,v_5]=\; w_2 & [v_3, v_5] = \; w_1

 \end{array}
 \end{equation}

  \begin{equation}\label{m_11}\tag{$\mg_{11}$}
 \begin{array}{lcl}
 [v_1,v_5]=\; w_2  & [v_3,v_6]=\; w_1

 \end{array}
 \end{equation}

 \begin{equation}\label{m_12}\tag{$\mg_{12}$}
 \begin{array}{lcl}

 [v_1,v_3]=\; w_3  & [v_1,v_5]=\; w_2

 \end{array}
 \end{equation}

 \begin{equation}\label{m_14}\tag{$\mg_{14}$}
 \begin{array}{lcl}

 [v_1,v_3]=\; w_3

 \end{array}
 \end{equation}

From Family $4$ and Table $5$:
 \begin{equation}\label{m_3}\tag{$\mg_{3}$}
 \begin{array}{lcl}

[v_5,v_3]=\; w_1 & [v_1,v_5]=\; w_2 & [v_3, v_1] = \; w_3

 \end{array}
 \end{equation}

\vspace{.3cm}

Note that according to \cite{GT} the corresponding algebras are pairwise non-isomorphic, meaning that for example, $\mathcal{U}_1 + \mg_{11}$ and  $\mathcal{U}_1 + \mg_{14}$ are non-isomorphic.}
\end{notation}

\noindent {\bf Anosov Lie algebras of type $(6, 3)$}.
 Let $\ngo$ be an Anosov Lie algebra of type
$(6,3)$ with no abelian factor, and let $\tau, \tau_1, \tau_2,$ be as in Proposition \ref{deco}.
As in the previous cases, we will denote by $X_i$ and $Y_j$ the eigenvectors of $\tau_1$ and $\tau_2$ respectively with corresponding
eigenvalues $\l_i$ and $\mu_j.$ Note that, by
using Lemma \ref{m1}, it is easy to see that the splitting of $\tau_1$ is $[6]$ or $[3;3] $ (see Definition \ref{splitting}).
We are going to study now each one of these cases separately.

\vspace{0.3cm}

 \noindent {\bf Case i:} The splitting of $\tau_1$ is $[6]$ or equivalently, the characteristic polynomial of  $\tau_1$ is irreducible over $\ZZ$. We note that, in particular, this implies that $\l_i \ne
 \l_j$ for all $i\ne j$. From this, it can be shown that one may assume that
\begin{equation}\label{631}
\begin{array}{lll} [X_1,X_2]=Y_1, &\quad [X_3,X_4]=Y_2,
&\quad[X_5,X_6]=Y_3.
\end{array}
\end{equation}

In fact, if this is not the case we can reorder the basis so that
\begin{equation}\label{632}
\begin{array}{lll}
 [X_1,X_2]=Y_1, &\quad& [X_1,X_3]=Y_2,
\end{array}
\end{equation}
and in this situation one would have to consider three possibilities for
$Y_3,$
\begin{equation}\label{633}
\begin{array}{lll}
 \mbox{\bf (I)}\;[X_1,X_4]=Y_3, & \quad \mbox{\bf (II)}\;[X_2,X_j]=Y_3,
&\quad \mbox{\bf (III)}\;[X_4,X_5]=Y_3.
\end{array}
\end{equation}

Note that each one of these situations represents a few others that
are totally equivalent to the one considered.

Case {\bf (III)} is the simplest one because directly  from (\ref{632}) it follows that  $\l_1^2\l_2\l_3\l_4\l_5=1$ and therefore we obtain the contradiction $\l_1=\l_6.$

Cases {\bf (I)} and {\bf (II)} can be done in a very similar way. That is, by considering the possibilities for the other brackets among the $X_i$ we get either a contradiction or (\ref{631}). Therefore, we can assume (\ref{631}):
$$\begin{array}{lll}
[X_1,X_2]=Y_1, &\quad [X_3,X_4]=Y_2, &\quad[X_5,X_6]=Y_3.
\end{array}
$$
If there are  no other nontrivial Lie brackets but these, then $\ngo
\simeq \hg_3\oplus\hg_3\oplus\hg_3 \simeq \mathcal{U}_1$ which is known to be Anosov (see
\cite{lw2}). If there are more nontrivial Lie brackets, without  any loss
of generality, we can assume that $[X_3,X_2]=c\;Y_3$ and moreover,
using that there is no abelian factor, one can  see that $\ngo_{\CC}$
is isomorphic to $\ngo'(a,b,c),$ given by
\begin{equation}\label{639}
\begin{array}{lll}
 [X_1,X_2]=Y_1, & [X_3,X_4]=Y_2 & [X_5,X_6]=Y_3,\\

[X_5,X_4]=a\;Y_1, & [X_1,X_6]=b\;Y_2 & [X_3,X_2]=c\;Y_3,\\
\end{array}
\end{equation}
for some $a,b,c \in \CC$.
It is easy to see that we can not add more nontrivial Lie brackets by using that $\tau$ is hyperbolic and $\ds \prod_{i = 1}^{6} \lambda_i = 1$.
To know which non-isomorphic Lie algebras we obtain from (\ref{639}), we start by noting that if $abc=0,$ from \cite{GT} we get that we have only two non-isomorphic Lie algebras:
 $$\ngo'(a,b,0) \simeq \ngo'(1,1,0) \simeq \mathcal{U}_1+\mg_{11},  \qquad \text{and} \qquad \ngo'(0,0,c)\simeq \ngo'(0,0,1) \simeq \mathcal{U}_1+\mg_{14},$$
   (see Notation \ref{Unotation}).

On the other hand, if $abc \ne 0,$  by changing the basis to
$$
\beta'=\{X_1,aX_2,\tfrac{1}{a}X_3,\tfrac{1}{c}X_4,
cX_5,X_6,aY_1,\tfrac{1}{ac}Y_2,cY_3\},
$$
we have that
$\ngo'(a,b,c) \simeq \ngo'(1,abc,1),$ and then (if $abc \ne 0$)
 $$\ngo'(a,b,c) \simeq \ngo'(s,s,s) \simeq \mathcal{U}_1+ s\; \mathcal{U}_2  \simeq s\;\mathcal{U}_2+ \mathcal{U}_3,$$
 \noindent where $s^3=abc \ne 0.$
Hence, from \cite{GT} we get the following non-isomorphic Lie algebras
\begin{equation}\label{su2u3}
\begin{array}{@{\bullet\;\;}ll}
 s\;\mathcal{U}_2+ \mathcal{U}_3,\;\text{ corresponding to }s^3 \ne \pm 1 &\text{ (Family 2)}\\
 \mathcal{U}_2+ \mathcal{U}_3,\;   \text{ corresponding to }s^3 = 1 &\text{ (Family 4)}\\
-\mathcal{U}_2+ \mathcal{U}_3,\; \text{ corresponding to }s^3 =-1 &\text{ (Family 5)}
\end{array}
\end{equation}

Here we are referring to the families from \cite[Section 4]{GT}.

We will consider now the ones corresponding to $s\in \mathbb{Q}$ which includes, in particular, the algebras  $\mathcal{U}_2+ \mathcal{U}_3$ and $-\mathcal{U}_2+ \mathcal{U}_3$.

\begin{proposition}\label{sU_2+U_3anosov}
If $s \in \QQ$, then a Lie algebra $ s\;\mathcal{U}_2+ \mathcal{U}_3$ is Anosov.
\end{proposition}

\proof   Let $s =\tfrac{p}{q} $ where $p, q \in \ZZ$ and $\ngo = s \;\mathcal{U}_2+ \mathcal{U}_3$ .
 Then $\ngo  \simeq p\,\mathcal{U}_2+q\,\mathcal{U}_3.$
 We recall that $p\,\mathcal{U}_2+q\,\mathcal{U}_3$ is a 2-step nilpotent Lie algebra in which the nonzero Lie brackets on its basis vectors
   $\{X_i, Y_j: 1\leq i \leq 6, 1 \leq j \leq 3 \}$  are  given by

   \begin{equation}
 \begin{array}{lcl}

[X_5,X_4]=p\;Y _1 & [X_1,X_6]=p\; Y_2 & [X_3, X_2]=p\; Y_3 \\

 [X_3,X_6]=q\; Y_1 & [X_5,X_2]=q\;Y_2 & [X_1,X_4]=q\;Y_3.
 \end{array}
 \end{equation}

  Let $\l_1, \l_2, \l_3$ be the roots of the polynomial $x^3 - 3x + 1$.  Note that we used these algebraic units in the proof of  Proposition \ref{(3,3,3)}.
Let $\sigma$  be the automorphism of  $\ngo$ whose matrix with respect to a basis $\{X_1,\ldots, X_6, Y_1, Y_2, Y_3 \}$ is a
diagonal matrix $D( \l_3,\l_3,\l_1,\l_1,\l_2,\l_2, \l_1\l_2,\l_2\l_3,\l_1\l_3)$.

 Let
 \begin{equation}\label{xyz63}
\begin{array}{l}
\mathcal{X}_l= \l_3^{l-1}\, X_1 + \l_1^{l-1}\, X_3+
\l_2^{l-1}\,X_5 \qquad l=1,2,3,\\
\mathcal{Y}_k= \l_3^{k-1}\, X_2 + \l_1^{k-1}\, X_4+
\l_2^{k-1}\,X_6 \qquad k=1,2,3,\\
 \mathcal{Z}_r= \l_1^r \,Y_1+\l_2^r \,Y_2+\l_3^r\, Y_3 \qquad r=1,2,3.
\end{array}\end{equation}
Hence,
   \begin{equation}\label{basis63}
\beta' = \{
\mathcal{X}_l, \mathcal{Y}_k, \mathcal{Z}_r \,|\; l,k,r=1,2,3 \}.
\end{equation}
 is a basis  of $\ngo$.  Using the properties of $\l_1, \l_2, \l_3$ it can be checked that
$[\sigma]_{\beta'} \in \Gl(9,\ZZ).$
  For example, as
$\l_1\l_2=\l_1-1$ (see (\ref{lambdas333})) we get that $\sigma(\mathcal{Z}_1) =  \mathcal{Z}_2 - \mathcal{Z}_1$. 

To see that $\beta'$ is a $\ZZ$-basis,  we note the following:

$$\begin{array}{l}
[\mathcal{X}_i,\mathcal{X}_j]= [\mathcal{Y}_i,\mathcal{Y}_j]= [\mathcal{Z}_i,\mathcal{Z}_j]=0 \quad \forall\,\,  1\le i,j \le 3,\\

[\mathcal{X}_1,\mathcal{Y}_1]= (p+q)\mathcal{Z}_1, \qquad [\mathcal{X}_2,\mathcal{Y}_1]=-2p\mathcal{Z}_1+q\mathcal{Z}_2 +p \mathcal{Z}_3,\\

 [\mathcal{X}_3,\mathcal{Y}_1]=2p\mathcal{Z}_1+p\mathcal{Z}_2 +q \mathcal{Z}_3.

\end{array}$$
\hfill$\square$

\vspace{.3cm}
\begin{remark}\label{infinite}
{\rm  Using the above proposition we get that $\{s \;\mathcal{U}_2+ \mathcal{U}_3: s \in \mathbb Q\}$ is an infinite family of non-isomorphic  Anosov Lie algebras of type $(6, 3)$. 
We know that there are, up to isomorphism, only finitely many (real) 
Anosov Lie algebras up to dimension 8 (see \cite[Table 3]{lw2}) and hence this is the smallest dimension in which there are infinitely many non-isomorphic Anosov Lie algebras.
  We also note that since  $\{s \;\mathcal{U}_2+ \mathcal{U}_3: s \not \in \QQ, s^3 \ne \pm 1\}$ is an uncountable family of non-isomorphic Lie algebras, not all of them can be Anosov. Moreover, if $s^3 \notin \QQ$ it is not even clear if $s \;\mathcal{U}_2+ \mathcal{U}_3$ admits a rational form. }
\end{remark}
\vspace{.3cm}

 {\bf Case ii:} The splitting of $\tau_1$ is $[3,3] $ (see Definition \ref{splitting}).  Let $\l_1, \ldots, \l_6$ denote the eigenvalues  of $\tau_1$ such that $\l_1, \l_2, \l_3$ are conjugates over $\QQ$ and $\l_4, \l_5, \l_6$ are conjugates
over $\QQ$.  We will follow  Notation \ref{eigenvectors}.

 In this case   one may either have
\begin{enumerate}
\item[{\bf (a)}] $[X_i,X_j]=Y_k$ for some $1\le i,j \le 3$ (or equivalently $4\le i,j \le 6$) or
\item[{\bf (b)}] $[X_i,X_j]=0$ for all $1\le i,j \le 3$ and $[X_k, X_l]= 0$ for all $4\le k, l \le 6$.
\end{enumerate}
In {\bf (a)}, with no loss of generality we may assume $[X_1,X_2]=Y_1.$  At
the eigenvalue level, this means that  $\lambda_1 \lambda_2=\l_3^{-1},\;
\lambda_1 \lambda_3=\l_2^{-1}$ and $\lambda_2 \lambda_3=\l_1^{-1}$ must be the
eigenvalues of $\tau_2$ ($\tau$ on $[\ngo, \ngo]$).
On the other hand, the absence of abelian factor implies that
$$
\begin{array}{lcl}
[X_4,X_j]\ne 0 &  \mbox{ for some } & 1\le j \le 6.
\end{array}$$

Hence $\l_4 \l_j = \l_i^{-1}$ for some $j$ and $i$,  $1\le j \le 6, \, 1 \le i \le 3$.
It is not hard to see that from here we can either have
$$\{ \l_4,\l_5,\l_6\} = \{ \l_1^{-2},\l_2^{-2},\l_3^{-2}\}
\qquad \mbox{or} \qquad \{ \l_4,\l_5,\l_6\} = \{\l_1,\l_2,\l_3\}. $$

If  {\bf $\{ \l_4,\l_5,\l_6\} = \{ \l_1^{-2},\l_2^{-2},\l_3^{-2}\}$},  we will prove that $\ngo$ is Anosov.  In this case, we may
 rearrange the basis so that  $\l_4 = \l_1^{-2},  \l_5 = \l_2^{-2} $ and $\l_6 = \l_3^{-2}$.
 Hence  we can assume that the  nonzero  Lie brackets in $\ngo$ are given by
\begin{equation}\label{6312}
\begin{array}{lll}
[X_1,X_2]=Y_1, & [X_2,X_3]=Y_2 & [X_1,X_3]=Y_3,\\

[X_3,X_6]=c\;Y_1, & [X_1,X_4]=a\;Y_2, & [X_2,X_5]=b\;Y_3  \\
\end{array}
\end{equation}

for some $a,b,c \in \CC$.  Here we need to use the special properties of $\l_i$'s like $|\l_i| \ne 1$, $\l_i$'s are all distinct  for $1 \le i \le 3$ etc.

  Let $\widetilde{\ngo}(a,b,c)$  denote the Lie algebra defined by (\ref{6312}).
  Due to our assumption of no abelian
factor we have that $abc \ne 0$ and then  $$\ngo \simeq \widetilde{\ngo}(a,b,c)
\simeq \widetilde{\ngo}(1,1,1)\simeq \mathcal{U}_1 + \mg_{10}$$ (see Notation \ref{Unotation} and (\ref{m_10})).

\begin{proposition}
A Lie algebra $ \widetilde{\ngo}(1,1,1)$ defined by (\ref{6312}) with $a= b= c= 1$ is Anosov.
\end{proposition}

\proof We will use very similar arguments as used in the case of type $(3, 3, 3)$ (proof of Proposition \ref{(3,3,3)}).  Let $\l_1, \l_2, \l_3$ be the
 roots of $X^3 - 3 X + 1.$ Let $\sigma$  be the automorphism of  $\widetilde{\ngo}(1,1,1)$ whose matrix with respect to  a basis $\{X_1,\ldots, X_6, Y_1, Y_2, Y_3 \}$ is a diagonal matrix $$D( \l_1,\l_2,\l_3,\l_1^{-2},\l_2^{-2},\l_3^{-2}, \l_3^{-1}, \l_1^{-1},\l_2^{-1}).$$

Let

\begin{align*}
\mathcal{X}_i &=  \l_1^{i-1} \, X_1+\l_2^{i-1} \, X_2+\l_3^{i-1} \, X_3, \\
\mathcal{X'}_j &= \l_3^{1-j}(\l_2-\l_1)X_6+ \l_1^{1-j}(\l_3-\l_2)X_4+\l_2^{1 - j}(\l_3-\l_1)X_5,  \\
\mathcal{Y}_k &= \l_3^{1-k}(\l_2-\l_1)Y_1+ \l_1^{1-k}(\l_3-\l_2)Y_2+\l_2^{1-k}(\l_3-\l_1)Y_3.
\end{align*}

It can be shown that $\tilde{\beta}=\{ \mathcal{X}_{i}, \mathcal{X'}_j,
\mathcal{Y}_k \,|\; i,j,k =1,2,3  \}$ is a $\ZZ$-basis of $\widetilde{\ngo}(1,1,1)$ such that $[\sigma]_{\widetilde{\beta}} \in \Gl(9,\ZZ)$ .
Note that due to the similarities with the $(3,3,3)$-case one only needs to check brackets that involve some $\mathcal{X}'_j$ (see \cite{mw}). \hfill$\square$

\vspace{0.3cm}

On the other hand,   if $\{\l_4, \l_5, \l_6 \} = \{\l_1, \l_2, \l_3\}$, then by rearranging the basis vectors we assume that $\l_4 = \l_1, \l_5 = \l_2, \l_6 = \l_3$ and
 the nonzero Lie brackets on $\ngo$ are given by
$$\begin{array}{lll}
[X_1,X_2]=a\;Y_1, & [X_2,X_3]=b\;Y_2 & [X_1,X_3]=c\;Y_3,\\

[X_4,X_5]=d\;Y_1, & [X_5,X_6]=e\;Y_2 & [X_4,X_6]=f\;Y_3,\\

[X_1,X_5]=g\;Y_1, & [X_2,X_6]=h\;Y_2 & [X_4,X_3]=i\;Y_3,\\

[X_4,X_2]=j\;Y_1, & [X_5,X_3]=k\;Y_2 & [X_1,X_6]=l\;Y_3,\\
\end{array}
$$

for some $a,b,c,d,e,f,g,h,i,j,k,l \in \CC.$
Note that the Pfaffian form of these Lie algebras (see \cite{lw1}) is given by
\begin{equation}
 Pf(x,y,z)=-xyz(afk-aei+bdl-bgf+cje-cdh+ihg-ljk).
 \end{equation}

Calculating the Pfaffian forms of the algebras listed in the classification given in \cite[Section 4]{GT},
one can see that $\ngo$ should be isomorphic to one of the following (see Notation \ref{Unotation}):
\begin{itemize}
\item $s\,\mathcal{U}_2+\mathcal{U}_3$,
\item $\mathcal{U}_2+\mathcal{U}_3+\mg_3$,
\item $\mathcal{U}_1+\mg_i$ for $i=6,7,10,11,12,14,15,$
\end{itemize}
\noindent where $\mg_{15}=0.$
\begin{remark}
  Note that the algebras in  \cite[Section 4, Family 7]{GT} with $0$ Pfaffian form  have always an abelian factor and therefore are not included in our list.
\end{remark}
The Lie algebras $s\,\mathcal{U}_2+\mathcal{U}_3$ have been considered in case {\bf i}  (see (\ref{su2u3}), Proposition \ref{sU_2+U_3anosov}).
   Also  $\mathcal{U}_1+\mg_i$ for $i=10,11,14,15$ have been considered in case {\bf i.}
   We don\rq{}t  know if  $\mathcal{U}_1+\mg_i$  is Anosov for $i=6,7,12$. However we prove the following:

\begin{proposition}
$\mathcal{U}_2+\mathcal{U}_3+\mg_3$ is an Anosov Lie algebra.
\end{proposition}

\proof Let $\ngo = \mathcal{U}_2+\mathcal{U}_3+\mg_3$. We recall that the nonzero brackets in $\ngo$
 on the basis vectors are given by (see Notation \ref{Unotation}):
$$\begin{array}{lll}
[X_5,X_4]=Y_1, & [X_1,X_6]=Y_2 & [X_3,X_2]=Y_3,\\

[X_3,X_6]=Y_1, & [X_5,X_2]=Y_2 & [X_1,X_4]=Y_3,\\

[X_5,X_3]=Y_1, & [X_1,X_5]=Y_2 & [X_3,X_1]=Y_3.
\end{array}$$

 Once again we will take $\l_1, \l_2, \l_3$ to be the roots of  $X^3 - 3X +1$. Let $\sigma$ be the automorphism of $\ngo$ whose matrix  with respect
 to the basis $\{X_i, Y_j: 1\le i \le 6, 1\le j \le 3 \}$ is a diagonal matrix
 $D( \l_3,\l_3,\l_1,\l_1,\l_2,\l_2, \l_1\l_2,\l_2\l_3,\l_1\l_3)$.
  Then $\sigma$ is a hyperbolic automorphism because of our choice of $\l_1, \l_2, \l_3$.

Let
 \begin{equation}
\begin{array}{l}
\mathcal{X}_l= \l_3^{l-1}\, X_1 + \l_1^{l-1}\, X_3+
\l_2^{l-1}\,X_5 \qquad l=1,2,3,\\
\mathcal{Y}_k= \l_3^{k-1}\, X_2 + \l_1^{k-1}\, X_4+
\l_2^{k-1}\,X_6 \qquad k=1,2,3,\\
 \mathcal{Z}_r=\l_3^{r-1}Y_1+ \l_1^{r-1}Y_2+ \l_2^{r-1}Y_3 \qquad r=1,2,3.
\end{array}\end{equation}

We note that a  similar construction of  $\mathcal{X}_l$'s, $\mathcal{Y}_k$'s and $\sigma$  have been used in the proof of Proposition \ref{sU_2+U_3anosov}.
Let  $\beta'= \{\mathcal{X}_l, \mathcal{Y}_k, \mathcal{Z}_r \,|\; l,k,r=1,2,3 \}$.  It can  seen  that $[\sigma]_{\beta'} \in \Gl(9,\mathbb{Z})$
 To show that $\beta'$ is a $\mathbb{Z}$-basis of $\ngo$, we note that
 $$ \lambda_1=\xi + \xi^8, \quad \lambda_2=\xi^2 + \xi^7,\quad \lambda_3=\xi^4 + \xi^5,$$
 where $\xi=e^{2i\pi/9}$ a ninth root of unity.

 Hence  if  $P(\l)=2\l^2+\l-4$,
  $$ \l_2-\l_3=P(\l_1),\quad\l_3-\l_1=P(\l_2),\quad \l_1-\l_2=P(\l_3).$$
Moreover,   $P(\l_i)\l_i=\l_i^2+2\l_i-2$ for $1\le i \le 3$.
Bearing all this in mind, straightforward calculations show that
$$\begin{array}{ll}
 [\mathcal{Y}_j,\mathcal{Y}_k]=[\mathcal{Z}_j,\mathcal{Z}_k]=0 \quad 1\le j,k \le 3,&[\mathcal{X}_1,\mathcal{Y}_1]=2\mathcal{Z}_1,  \\

 [\mathcal{X}_1,\mathcal{X}_2]=-4\mathcal{Z}_1+\mathcal{Z}_2+2\mathcal{Z}_3,& [\mathcal{X}_1,\mathcal{Y}_2]=[\mathcal{X}_2,\mathcal{Y}_1]=-\mathcal{Z}_2,\\

 [\mathcal{X}_1,\mathcal{X}_3]= 2\mathcal{Z}_1-2\mathcal{Z}_2-\mathcal{Z}_3,& [\mathcal{X}_1,\mathcal{Y}_3]=[\mathcal{X}_3,\mathcal{Y}_1]=6\mathcal{Z}_1-\mathcal{Z}_2.
\end{array}
$$
Hence   $\mathcal{U}_2+\mathcal{U}_3+\mg_3$ is an Anosov Lie algebra. \hfill$\square$

\vspace{.3cm}

To conclude, let us study {Case \bf (b)}. We recall that in this case,  our assumption is
$$[X_i,X_j]=0=[X_k, X_l]  \,\, \mbox{ for }\,\,1 \le i,j
\le 3,\;\; 4 \le k, l \le 6.$$
We  will show that $\ngo$ is  isomorphic to one of the Lie algebras we came across  in Case (a).  In other words,
 we don't get any new  Lie algebras in this case.

 We note that we can assume
\begin{equation}\begin{array}{lcl}\label{u1}
[X_1,X_4]=Y_1, & [X_2,X_5]=Y_2 & [X_3,X_6]=Y_3,
\end{array}
\end{equation}
In fact, if this is not the case, rearranging a basis we can only assure that
$$\begin{array}{lll}
[X_1,X_4]= Y_1, & [X_i,X_j]=Y_2 & [X_k,X_l]=Y_3,
\end{array}
$$
$i,k \le 3$ $j,l \ge 4$ and $i$ or $k =1.$ Let us say that $[X_1,X_5]=Y_2$ and then it is not hard to check that $k \ne 1$ so we may assume that  $k = 2.$ Also, since $l\ne 6$ ($l=6$ implies the contradiction $\l_1=\l_3$) , say $l=5$ ($l = 4$ is similar). Then we get  at least the following nonzero  Lie brackets:
$$\begin{array}{lll}
[X_1,X_4]= Y_1, & [X_1,X_5]=Y_2 & [X_2,X_5]=Y_3.
\end{array}
$$
Therefore $\l_1^2\l_2\l_4\l_5^2=1$ or equivalently, $\l_1\l_5=\l_3\l_6$ and from
this, since $\ngo$ has no abelian factor, by considering all possibilities for $[X_r,X_6]=Y_k$ and $[X_3,X_s]=Y_l$ it is not hard to check that we should have $[X_3,X_6]=Y_2,$ as desired, since any other possibility leads to a contradiction.

If $\ngo$ has only the nonzero Lie brackets given in (\ref{u1}) then $\ngo \simeq \mathcal{U}_1$ which is Anosov.
If there are more nonzero Lie brackets, with no loss of
generality, we may assume that $[X_3,X_5]= a\;Y_1,$ and hence we
have that
$$ \mu_1=\l_1\l_4=\l_3\l_5.$$ Note that since
 $\l_i$'s are distinct for $ 1 \le i \le 3$ and $\l_j$\rq{}s are distinct for $4 \le j \le 6$,
 if $[X_i,X_j]= Y_k=[X_l,X_r]$
then $i \ne l$ and $j \ne r$. From this, it can be seen that
 the nonzero  Lie brackets  in  $\ngo$ are given by
\begin{equation}\label{6314}
\begin{array}{lll}
[X_1,X_4]= Y_1, & [X_2,X_5]=Y_2 & [X_3,X_6]=Y_3,\\ &&\\

[X_3,X_5]= a\;Y_1, & [X_1,X_6]=b\;Y_2 & [X_2,X_4]=c\;Y_3,
\end{array}
\end{equation}

for some constants $a,b,c \; \in \CC.$ By reordering the
basis
$$\beta'=\{X_1,X_4,X_2,X_5,X_3,X_6,Y_1,Y_2,Y_3\}$$
one can see that this algebra is isomorphic to $\ngo(a, b, c)$ given in (\ref{639}) (from Case (a))  that has already been studied.

\vspace{.3cm}

Finally, if $\ngo$ is Anosov and has an abelian factor,  and if   $\ngo=\widetilde{\ngo}\oplus\mg,$ where
$\mg$ is a maximal abelian  factor of $\ngo$, according to \cite[Theorem 3.1]{lw1} one has that $\dim \mg\geq 2.$ Moreover, $\widetilde{\ngo}$ is also an Anosov Lie algebra and $\dim \widetilde{\ngo}\le 7$. Since there are  no $7$ dimensional  non-abelian Anosov Lie algebras,  we then have that $\ngo$ is isomorphic to one of the following  (see \cite{lw2}):
\begin{itemize}
\item $\RR^9$
\item $\hg_3 \oplus \hg_3 \oplus \RR^3$
\item $\fg_3 \oplus \RR^3, $
\end{itemize}
\noindent where $\hg_3$ is the $3$-dimensional Heisenberg algebra,
and $\fg_3$  is the free $2$-step nilpotent Lie algebra on 3
generators.

Summarizing the results of this section we have the following Theorem.

\begin{theorem} \label{teo}If $\ngo$ is a 9-dimensional Anosov Lie algebra then
\begin{description}
  \item[i] if it has no abelian factor, then one of the following holds
  \begin{itemize}
    \item it is a type $(3,3,3)$ nilpotent Lie algebra isomorphic to $\ngo(1,1,1)$ given in (\ref{n333}) or
    \item it is a two step nilpotent Lie algebra of type $(6,3)$ and its Pfaffian
form is projectively equivalent to $xyz$ or $0$. Moreover it
 is isomorphic to one of the following non-isomorphic Lie
 algebras: $\mathcal{U}_1,$ $\mathcal{U}_1+\mg_i$ for
$i=6,7,10,11,12,14,$ $s\,\mathcal{U}_2+\mathcal{U}_3$ $s \in \CC$ or $\mathcal{U}_2+\mathcal{U}_3+\mg_3.$
  \end{itemize}
  \item[ii] If $\ngo$ has an abelian factor then it is isomorphic to $\RR^9,$
 $\hg_3 \oplus \hg_3 \oplus \RR^3$ or $\fg_3 \oplus \RR^3.$
\end{description}

\end{theorem}

\noindent{\bf Acknowledgements. } The authors would like to thank J. Lauret for his constant help and C. Gordon for her hospitality and encouragement during this project.  Meera  Mainkar was supported by the Central Michigan University ORSP Early Career Investigator (ECI) grant \#C61940.  Cynthia Will was supported by  CONICET and grants from FONCyT and SeCyT.

\end{document}